\documentclass{article}
\usepackage{amssymb}
\usepackage{amsmath}
\newtheorem{theorem}{Theorem}

\newtheorem{lemma}{Lemma}
\newtheorem{corollary}{Corollary}

\title{A collection of sharp  dilation invariant inequalities for differentiable functions}
\author{
Vladimir Maz'ya\thanks{The first author was partially  
supported by the USA National Science Foundation grant DMS 
0500029 and by the UK Engineering and Physical Sciences Research Council  grant EP/F005563/1.}
\\
\emph{\small Department of Mathematical Sciences,  The University of Liverpool, Liverpool L69 7ZL}\\\emph{\small Department of Mathematics, Link\"oping University, Link\"oping SE-58183, Sweden}
\\
\\*[10pt]
Tatyana Shaposhnikova\\
\emph{\small Department of Mathematics, Link\"oping University, Link\"oping SE-58183, Sweden}\\
}
\date{}

\begin{document}
\maketitle
\centerline
{\it In memory of S. L. Sobolev}

\bigskip

\medskip

\begingroup

\narrower\noindent
{\it Abstract.}  We find best constants in several dilation invariant integral inequalities involving derivatives of functions. Some of these inequalities are new and some were known without best constants. The contents: 1. Estimate for a quadratic form of the  gradient, 2. Weighted G{\aa}rding inequality for the biharmonic operator, 3. Dilation invariant Hardy's inequalities  with  remainder term, 4. Generalized Hardy-Sobolev inequality with sharp constant, 5. Hardy's inequality with sharp Sobolev remainder term.
\endgroup



\bigskip

\section*{Introduction}

The present article consists of five independent sections dealing with various dilation invariant integral  inequalities with optimal constants. We briefly describe the contents, starting with Section 1. 

\smallskip

Let us recall  the Gagliardo-Nirenberg inequality
\begin{equation}\label{a1}
\|v\|_{L^2(\Bbb{R}^2)} \leq C \, \|\nabla v\|_{L^1(\Bbb{R}^2)}, \quad v\in C_0^\infty(\Bbb{R}^2)
\end{equation}
\cite{[G]}, \cite{[N]}. (The best constant $C= (2\sqrt{\pi})^{-1}$ was found in \cite{[FF]} and  \cite{[M1]}, see also \cite{[M5]}, Sect. 1.4.2). Setting 
 $v= |\nabla u|$, we observe that the Dirichlet integral of $u$ admits the estimate
\begin{equation}\label{x1}
\int_{\Bbb{R}^2}|\nabla u |^2 dx  \leq C\, \Bigl( \int_{\Bbb{R}^2} |\nabla_2 u|\, dx\Bigr)^2,
\end{equation}
where
$$|\nabla_2 u|^2 = |u_{x_1x_1}|^2 + 2 |u_{x_1x_2}|^2 + |u_{x_2x_2}|^2.$$
One can see that it is impossible to improve   (\ref{x1}), replacing $|\nabla_2 u|$ in the right-hand side by $|\Delta u|$. Indeed, it suffices to put a sequence of mollifications of the function $x\to \eta(x)\log |x|$, where $\eta\in C_0^\infty(\Bbb{R}^2)$, $\eta(0) \neq 0$, into the estimate in question in order to check its failure.

\smallskip

However, we show that the estimate of the same nature
$$\Bigl |\int_{\Bbb{R}^2 } \sum_{i,j =1}^2 a_{i,j}\, u_{x_i}\, {\overline{u}}_{x_j}\, dx \Bigr |\leq C\, \Bigl( \int_{\Bbb{R}^2} |\Delta u|\, dx\Bigr)^2,$$
where $a_{i,j} = const$ and $u$ is an arbitrary complex-valued function in $C_0^\infty(\Bbb{R}^2)$,  may hold if and only if 
$a_{11} + a_{22} =0$. We also find the best constant $C$ in the last inequality. 
This is a particular case of Theorem \ref{th5} proved in Section 1.

\smallskip

In Section 2 we establish a new weighted G{\aa}rding type inequality
\begin{equation}\label{t1}
\int_{\Bbb{R}^2 } |\nabla_2 u|^2 \log(e^2|x|)^{-1}dx \leq {\rm Re} \int_{\Bbb{R}^2 }\Delta^2u\cdot \overline{u} \, \log |x|^{-1} dx
\end{equation}
for all $u\in C_0^\infty(\Bbb{R}^2 \backslash \{0\})$. Estimates of such a kind proved to be useful in the study of boundary behavior of solutions to elliptic equations (see \cite{[M4]},   \cite{[MT]},  \cite{[M6]}, \cite{[M7]}, \cite{[E]}, \cite{[MM]}).

\smallskip

Before turning to the contents of the next section, we introduce some notation. By $\Bbb{R}^n_+$ we denote the half-space $\{x =(x_1, \ldots x_n)\in \Bbb{R}^n, x_n>0\}$. Also let $\Bbb{R}^{n-1} = \partial \Bbb{R}^n_+$. As usual, $C_0^\infty(\Bbb{R}^n_+)$ and $C_0^\infty(\overline{\Bbb{R}^n_+})$ stand for the spaces of infinitely differentiable functions with compact support in $\Bbb{R}^n_+$ and $\overline{\Bbb{R}^n_+}$, respectively.

\smallskip

In Section 3 we are concerned with the inequality
\begin{equation}\label{1}
\int_{\Bbb{R}^n_+} x_n\, |\nabla u|^2 dx \geq \Lambda\int_{\Bbb{R}^n_+}\frac{|u|^2}{(x_{n-1}^2 + x_n^2)^{1/2}} dx, \quad u\in C_0^\infty(\overline{\Bbb{R}^n_+}).
\end{equation}
It was obtained in  1972 by one of the authors   
 and proved to be useful in the study of the generic case of degeneration in the  oblique derivative problem for  second order elliptic differential operators  
 \cite{[M2]}. 

\smallskip

 By  substituting $u(x) = x_n^{-1/2} v(x)$ 
 into (\ref{1}),  one deduces with the same $\Lambda$ that
\begin{equation}\label{2}
\int_{\Bbb{R}^n_+}  |\nabla v|^2 dx \geq \frac{1}{4} \int_{\Bbb{R}^n_+} \frac{|v|^2 dx}{x_n^2} + \Lambda \int_{\Bbb{R}^n_+}\frac{|v|^2 dx}{x_n(x_{n-1}^2 + x_n^2)^{1/2} } 
\end{equation}
for all $v\in C_0^\infty(\Bbb{R}^n_+)$ (see \cite{[M5]},   Sect. 2.1.6).

Another inequality of a similar nature obtained in \cite{[M5]} is 
\begin{equation}\label{6x}
\int_{\Bbb{R}^n_+} |\nabla v|^2 dx \geq \frac{1}{4}\int_{\Bbb{R}^n_+} \frac{|v|^2}{x_n^2} dx+ C\,  \|x_n^\gamma \,  v\, \|^2_{L_q(\Bbb{R}^n_+)}.
\end{equation}
(This is a special  case of inequality (2.1.6/3) in \cite{[M5]}.)

\smallskip

 Without the second term in the right-hand sides of (\ref{2}) and (\ref{6x}), these inequalities reduce to the classical Hardy inequality with the sharp constant $1/4$ (see \cite{[Da]}).   An interesting feature of (\ref{2}) and (\ref{6x}) is their  dilation invariance. 
 
Variants, extensions, and  refinements of (\ref{2}) and (\ref{6x}), usually  called Hardy's inequalities with remainder term,   became the theme of many subsequent studies (\cite{[ACR]}, \cite{[Ad]}, \cite{[AGS]}, \cite{[BCT]}, \cite{[BFT1]},  \cite{[BFT2]}, \cite{[BFL]}, \cite{[BM]}, \cite{[BV]},  \cite{[CM]}, \cite{[DD]}, \cite{[DNY]}, \cite{[EL]}, \cite{[FMT1]}--\cite{[FMT3]}, \cite{[FT]}, \cite{[FTT]}, \cite{[FS]}, \cite{[GGM]},  \cite{[HHL]},   \cite{[TT]}, \cite{[TZ]}, \cite{[Ti1]}, \cite{[Ti2]}, \cite{[VZ]}, \cite{[YZ]} {\it et al}). 

\smallskip

In Theorem \ref{prop1}, proved in Section 3, we find a  condition on the function $q$ which is necessary and sufficient for  the inequality
\begin{equation}\label{2u}
\int_{\Bbb{R}^n_+}  |\nabla v|^2 dx - \frac{1}{4} \int_{\Bbb{R}^n_+} \frac{|v|^2 dx}{x_n^2} \geq  C \int_{\Bbb{R}^n_+}q \Bigl( \frac{x_n}{(x_{n-1}^2 + x_n^2)^{1/2}}\Bigr)
\frac{|v|^2 dx}{x_n\, (x_{n-1}^2 + x_n^2)^{1/2} }, 
\end{equation}
where $v$ is an arbitrary function in $C_0^\infty(\Bbb{R}^n_+)$. 
This condition  implies, in particular, that  the right-hand side of (\ref{2}) can be replaced by
$$C \int_{\Bbb{R}^n_+}\frac{|v|^2 dx}{x_n^2\Bigl(1-\log \displaystyle{\frac{x_n}{(x_{n-1}^2 + x_n^2)^{1/2}}}\Bigr)^2}.$$

\smallskip

The value $\Lambda =1/16$  in (\ref{1}) obtained in   \cite{[M2]} is not best possible.  Tidblom replaced it by $1/8$ in \cite{[Ti2]}. As a corollary of Theorem \ref{prop1}, we find an expression for the optimal value of $\Lambda$. 

\smallskip

Let the measure $\mu_b$ be defined by
 \begin{equation}\label{65j}
 \mu_b(K) = \int_K\frac{dx}{|x|^b}
 \end{equation}
for any compact set $K$ in $\Bbb{R}^n$. In 
 Section 4 we obtain the best constant in the inequality
 $$\|u\|_{{\cal L}_{\tau,q}(\mu_b)} \leq C \Bigl(\int_{\Bbb{R}^n}{| \nabla u(x)|^p}\frac{dx}{|x|^a} \Bigr)^{{1}/{p}},$$
where the left-hand side is the quasi-norm in the Lorentz space ${\cal L}_{\tau,q}(\mu_b)$, i.e.
$$\|u\|_{{\cal L}_{\tau,q}(\mu_b)} =  \Bigl( \int_0^\infty \bigl(\mu_b\{x: |u(x)|\geq t\} \bigr)^{q/\tau} d(t^q) \Bigr)^{1/q}.$$
As a particular case of this result we obtain the best constant in the
 Hardy-Sobolev inequality 
 \begin{equation}\label{NN}
\left(\int_{\Bbb{R}^n}{|u(x)|^q}\frac{dx}{|x|^b} \right)^{{1}/{q}}\leq {\cal C} \Bigl(\int_{\Bbb{R}^n}{| \nabla u(x)|^p}\frac{dx}{|x|^a} \Bigr)^{{1}/{p}},
\end{equation}
first proved by Il'in  in 1961 in \cite{[Il]} (Th. 1.4) without discussion of  the value of ${\cal C}$.
 Our result is  a direct consequence of the capacitary integral inequality from \cite{[M7]} combined with an isocapacitary inequality. For particular cases the best   constant ${\cal C}$ was found in \cite{[CW]} ($p=2$), in 
 \cite{[M3]}, Sect. 2 ($p =1$, $a=0$), in \cite{[GMGT]} ($p=2$, $n=3$, $a=0$), 
in \cite{[L]} ($p=2$, $n\geq 3$, $a=0$), and in \cite{[Na]} ($1<p<n$, $a=0$), where different methods were used.

\smallskip

The topic of the concluding Section 5 is the best constant $C$ in 
 the inequality (\ref{6x}), 
where $u\in C^\infty(\overline{\Bbb{R}^n_+})$, $u=0$ on $\Bbb{R}^{n-1}$.

\smallskip

Recently Tertikas and Tintarev \cite{[TT]} obtained (among other results) the existence of an optimizer in (\ref{6x}) in the case $\gamma =0$, $q= 2n/(n-2)$, $n\geq 4$. However, for these values of $\gamma$, $q$, and $n$ the best value of $C$ is unknown. In the case $n=3$, $\gamma =0$, $q=6$  Benguria, Frank, and Loss  proved the nonexistence of an optimizer and found the best value of $C$ by an ingenious argument \cite{[BFL]}. 

\smallskip

We note in Section 5  that a similar problem can be easily solved  for the special case $q= 2(n+1)/(n-1)$ and $\gamma =-1/(n+1)$.

\section{Estimate for a quadratic form of the  gradient}

\begin{theorem}\label{th5}
Let $n\geq 2$ and let $A= \|a_{i,j} \|_{i,j =1}^n$ be an arbitrary  matrix with constant complex entries. The inequality
\begin{equation}\label{x2}
\Bigl |\int_{\Bbb{R}^n }\langle A\nabla u, \nabla u\rangle _{\Bbb{C}^n}\, dx \Bigr | \leq C  \Bigl( \int_{\Bbb{R}^n }\bigl | (- \Delta)^{\frac{n+2}{4}} u\bigr |\, dx \Bigr)^2, 
\end{equation}
where $C$ is a positive constant, holds for all complex-valued $u\in C_0^\infty(\Bbb{R}^n)$ if and only if the trace of $A$ is equal to zero. The best value of  $C$  is given by
\begin{equation}\label{x3}
C = \frac{(4\pi)^{-n/2}}{\Gamma\bigl(\frac{n}{2} +1\bigr)} \max\limits_{\omega\in S^{n-1}} \Bigl |\sum_{1\leq i,j\leq n} a_{i,j} \, \omega_i\omega_j \Bigr |,
\end{equation}
where $S^{n-1}$ is the $(n-1)$-dimensional unit sphere in $\Bbb{R}^n$.
\end{theorem}

(The notation $(-\Delta)^s$ in (\ref{x2}) stands for an  integer or noninteger power of $-\Delta$.)

\medskip

\noindent
{\bf Proof.} By ${\cal F}$ we denote the unitary Fourier transform in $\Bbb{R}^n$ defined by
 \begin{equation}\label{ax}
 {\cal F} h(\xi) = (2\pi)^{-n/2}\int_{\Bbb{R}^n }  h(x)\, e^{-i\, x \cdot \xi} \, dx.
\end{equation}
We set $h= (-\Delta)^{(n+2)/4} u$ and write  (\ref{x2}) in the form
\begin{equation}\label{x4}
\Bigl |\int_{\Bbb{R}^n } |{\cal F} h(\xi)|^2\, \Bigl\langle A\frac{\xi}{|\xi|}, \frac{\xi}{|\xi|}\Bigr\rangle_{\Bbb{C}^n}\, \frac{d\xi}{|\xi|^n}\Bigr | \leq C \Bigl(\int_{\Bbb{R}^n } |h(x)|\, dx\Bigr)^2.
\end{equation}
The singular integral in the left-hand side exists in the sense of the Cauchy principal value,  since
$$\int_{S^{n-1}} \langle A\omega, \omega\rangle _{\Bbb{C}^n}\, ds_\omega = n^{-1} |S^{n-1}| \, {\rm Tr} \, A =0,$$
where ${\rm Tr}\,  A$ is the  trace of $A$ (see, for example,  \cite{[MP]}, Ch. 9, Sect. 1 or \cite{[SW]}, Theorem 4.7). Let
$$k(\xi) = |\xi|^{-n} \Bigl\langle A\frac{\xi}{|\xi|}, \frac{\xi}{|\xi|}\Bigr\rangle_{\Bbb{C}^n}.$$
The left-hand side in (\ref{x4}) equals
$$\Bigl |\int_{\Bbb{R}^n } {\cal F}^{-1} \Bigl(k(\xi) \bigl ({\cal F} h\bigr)(\xi)\Bigr) (x) \, \overline{h(x)}\, dx \Bigr | = (2\pi)^{-n/2} \Bigl |\int_{\Bbb{R}^n }\bigl ( ( {\cal F}^{-1} k)\ast h\bigr) (x) \, \overline{h(x)}\, dx\Bigr |
$$
with $\ast$ meaning the convolution. Thus, inequality (\ref{x4}) becomes
\begin{equation}\label{Na}
\Bigl |\int_{\Bbb{R}^n }\bigl ( ( {\cal F}^{-1} k)\ast h\bigr) (x) \, \overline{h(x)}\, dx\Bigr | \leq (2\pi)^{n/2}\, C
\Bigl(\int_{\Bbb{R}^n } |h(x)|\, dx\Bigr)^2.
\end{equation}

\smallskip

We note that for $\xi\in \Bbb{R}^n $, 
\begin{equation}\label{x4a}
k(\xi) = |\xi|^{-n-2} \Bigl( \sum_{j=1}^n a_{jj} \bigl( \xi^2_j - n^{-1} |\xi|^2\bigr) + \sum_{i, j=1\atop{ i\neq j}}^n a_{ij}\xi_i\xi_j\Bigr).
\end{equation}
Hence, for $n>2$
\begin{equation}\label{x4b}
k(\xi)=\frac{1}{n\, (n-2)} \sum_{i,j=1}^n a_{ij}\frac{\partial ^2}{\partial\xi_i\partial \xi_j} |\xi|^{2-n},
\end{equation}
and for $n=2$
\begin{equation}\label{x4c}
k(\xi)=\frac{1}{2} \sum_{i,j=1}^2 a_{ij}\frac{\partial ^2}{\partial\xi_i\partial \xi_j}\log |\xi|^{-1}.
\end{equation}
Applying ${\cal F}^{-1}$ to the identity
$$-\Delta_\xi \Bigl( \frac{\partial ^2}{\partial\xi_i\partial \xi_j} \frac{|\xi|^{2-n}}{|S^{n-1}|\, (n-2)}\Bigr ) = \frac{\partial ^2}{\partial\xi_i\partial \xi_j} \delta(\xi),$$
where $n>2$ and $\delta$ is the Dirac function, we obtain from (\ref{x4b})
\begin{equation}\label{x5}
\bigl( {\cal F}^{-1} k\bigr)(x)= \frac{-|S^{n-1}|}{n\, (2\pi)^{n/2}} \sum_{i,j=1}^n a_{ij}\frac{x_ix_j}{|x|^2}.
\end{equation}
Here $|S^{n-1}|$ stands for the $(n-1)$-dimensional measure of  $S^{n-1}$:
\begin{equation}\label{nn1}
|S^{n-1}| = \frac{2\pi^{\frac{n}{2}}}{\Gamma(\frac{n}{2})}.
\end{equation}
Hence
\begin{equation}\label{nn2}
\bigl( {\cal F}^{-1} k\bigr)(x)= \frac{-2^{-n/2}}{\Gamma(1+\frac{n}{2})}\sum_{i,j=1}^n a_{ij}\frac{x_ix_j}{|x|^2}.
\end{equation}
Similarly,  we deduce from (\ref{x4c}) that (\ref{nn2}) holds for $n=2$ as well. 
%\begin{equation}\label{nn2}
%\bigl( {\cal F}^{-1} k\bigr)(x) =\frac{-1}{2} \sum_{i,j=1}^2 a_{ij}\frac{x_ix_j}{|x|^2}.
%\end{equation}
Now, (\ref{x2}) with $C$ given by (\ref{x3}) follows from (\ref{nn2})  inserted into (\ref{Na}).

\smallskip

Next, we show the sharpness of $C$ given by (\ref{x3}). Let $\theta$ denote a point  on $S^{n-1}$ such that 
\begin{equation}\label{x7}
\bigl |\bigl( {\cal F}^{-1} k\bigr)(\theta)\bigr | = \max\limits_{\xi\in \Bbb{R}^n\backslash\{0\}} |{\cal F}^{-1} k(\xi)| .
\end{equation}
In order to obtain the required lower estimate for $C$, it suffices to set 
$$h(x)= \eta(|x|)\,  \delta_\theta\Bigl(\frac{x}{|x|}\Bigr),$$
where $\eta\in C_0^\infty[0, \infty)$, $\eta\geq 0$, and $\delta_\theta$ is the Dirac measure on $S^{n-1}$  concentrated at $\theta$, into the inequality (\ref{Na}). 
 (The legitimacy of this choice  of $h$ can be easily checked by approximation.) Then  estimate (\ref{Na}) becomes
$$
\Bigl |\int_0^\infty\!\int_0^\infty \bigl( {\cal F}^{-1} k\bigr)\Bigl(\frac{\rho -r}{|\rho -r|}\theta\Bigr) \eta(r)\, r^{n-1} \, \eta(\rho)\, \rho^{n-1} dr d\rho\Bigr| $$
\begin{equation}\label{x8}
\leq (2\pi)^{n/2}\, C \Bigl(\int_0^\infty\eta(\rho)\, \rho^{n-1} d\rho \Bigr)^2.
\end{equation} 
In view of (\ref{x5}) and (\ref{nn2}),
$$\bigl( {\cal F}^{-1} k\bigr)(\pm\theta) = \bigl( {\cal F}^{-1} k\bigr)(\theta)$$
which together with (\ref{x7})  enables one to write (\ref{x8}) in the form
$$\max\limits_{\xi\in \Bbb{R}^n\backslash\{0\}} |{\cal F}^{-1} k(\xi)|  \leq (2\pi)^{n/2}\, C.$$
By (\ref{x5}) and (\ref{nn2}) this can be written as
$$\frac{|S^{n-1}|}{n\, (2\pi)^{n/2}} \max\limits_{\omega\in S^{n-1}} \Bigl |\sum_{1\leq i,  j\leq n} a_{i j} \, \omega_i\omega_j \Bigr |
 \leq (2\pi)^{n/2}\, C.$$
The result follows from (\ref{nn1}). $\square$

\bigskip

\noindent
{\bf Remark 1} Let $P$ and $Q$ be functions,  positively homogeneous of degrees $2m$ and $m + n/2$ respectively, $m>-n/2$. We assume that the restrictions of $P$, $Q$, and $P|Q|^{-2}$ to $S^{n-1}$ belong to  $L^1(S^{n-1})$, By the same argument  as in Theorem \ref{th5},  
one  concludes that the condition
\begin{equation}\label{x9}
\int_{S^{n-1}} \frac{P(\omega)}{|Q(\omega)|^2} \, ds_\omega =0
\end{equation}
is necessary and sufficient for the inequality
\begin{equation}\label{x10}
\Bigl |\int_{\Bbb{R}^n }  P(D) u \cdot \overline{u}\, dx \Bigr |\leq C \Bigl( \int_{\Bbb{R}^n }\bigl | Q(D)\,  u\bigr |\, dx \Bigr)^2
\end{equation}
to hold for all $u\in C_0^\infty(\Bbb{R}^n )$. Moreover,  using the classical formula for the Fourier transform of a  positively homogeneous function of degree  $-n$ (see Theorem 4.11 in \cite{[SW]}\footnote{Note that the definition of the Fourier transform in \cite{[SW]} contains ${\rm exp}( - 2\pi i\, x\cdot \xi)$ unlike (\ref{ax}).}), one obtains that the best value of $C$ in (\ref{x10}) is given by
\begin{equation}\label{y10}
\sup\limits_{\omega\in S^{n-1}}\Bigl |\int_{S^{n-1}} \Bigl(\frac{i\pi}{2} {\rm sgn}(\theta\cdot \omega) + \log |\theta\cdot \omega|\Bigr)   \frac{P(\theta)}{|Q(\theta)|^2} \, ds_\theta \Bigr |.
\end{equation}
In particular, if  $P(\omega)/|Q(\omega)|^2$ is a spherical harmonic,
 the best value of $C$  in (\ref{x10}) is equal to
\begin{equation}\label{z10}
\frac{(4\pi)^{-n/2} \Gamma(m)}{\Gamma\bigl(\frac{n}{2} + m\bigr)} \max\limits_{\omega\in S^{n-1}}
\frac{|P(\omega)|}{|Q(\omega)|^2}, 
\end{equation}
which coincides with (\ref{x3}) for $m=1$, $P(\xi) = A\xi \cdot \xi$, and $Q(\xi) = |\xi |^{1+n/2}$.

\section{Weighted G{\aa}rding inequality for the biharmonic operator}

We start with an auxiliary Hardy type inequality.

\begin{lemma}\label{lem4}
Let $u\in C_0^\infty(\Bbb{R}^2)$. Then the sharp inequality
\begin{equation}\label{1m}
\Bigl |{\rm Re}\int_{\Bbb{R}^2}\bigl(x_1 u_{x_1} + x_2 u_{x_2}\bigr) \Delta\overline{ u} \, \frac{dx}{|x|^2}\Bigr | \leq \int_{\Bbb{R}^2} |\Delta u|^2 dx
\end{equation}
holds.
\end{lemma}

\noindent
{\bf Proof.} Let $(r, \varphi)$ denote polar coordinates in $\Bbb{R}^2$ and let 
$$u(r,\varphi) = \sum_{k=-\infty}^\infty u_k(r) e^{ik\varphi}.$$
Then (\ref{1m}) is equivalent to the sequence of inequalities
\begin{equation}\label{2m}
\Bigl |{\rm Re} \int_0^\infty \Bigl( v'' + \frac{1}{r}v' - \frac{k^2}{r^2} v\Bigr) \overline{v}'\, dr\Bigr | \leq 
 \int_0^\infty \Bigl | v'' + \frac{1}{r}v' - \frac{k^2}{r^2} v\Bigr | ^2 r\, dr, \,\,\, k= 0,1,2,\ldots
\end{equation}
where $v$ is an arbitrary function on $C_0^\infty([0,\infty))$.  Putting $t= \log r^{-1}$ and $w(t) = v(e^{-t})$, we write (\ref{2m}) in the form
$$
\Bigl |{\rm Re} \int_{\Bbb{R}^1}\bigl( w'' -k^2w\bigr)\, \overline{w}'\, e^{2t}\, dt\Bigr | \leq \int_{\Bbb{R}^1}\bigl | w'' -k^2 w \bigr |^2 \, e^{2t}\, dt$$
which is equivalent to the inequality
\begin{equation}\label{3m}
\Bigl |{\rm Re} \int_{\Bbb{R}^1} \bigl(g'' -2g' + (1-k^2)g\bigr) \bigl(\overline{g}' - \overline{g}\bigr) dt\Bigr | \leq  \int_{\Bbb{R}^1} \bigl |g'' -2g' + (1-k^2)g \bigr |^2 dt,
\end{equation}
where $g= e^tw$. Making use of the Fourier transform in $t$, we see that (\ref{3m}) holds if and only if for all $\lambda\in \Bbb{R}^1$ and $k=0,1,2 \ldots$
$$
\Bigl |{\rm Re}\, \bigl (-\lambda^2 +1 -k^2 -2i\lambda\bigr)\bigl(1-i\lambda\bigr)\Bigr | \leq \bigl(\lambda^2 -1 + k^2\bigr)^2 + 4\lambda^2,$$
which is the same as
$$\bigl |3x -1 +k^2 \bigr | \leq x^2 +2(k^2 +1)x + (k^2-1)^2$$
with $x=\lambda^2$. This elementary inequality  becomes equality if and only if $k=0$ and $x=0$. $\square$

\bigskip

\noindent
{\bf Remark 2} In spite of the simplicity of its proof, inequality (\ref{1m}) deserves some interest. Let us denote the integral over $\Bbb{R}^2$ in the left-hand side of (\ref{1m}) by $Q(u,u)$ and write (\ref{1m}) as
$$|\, {\rm Re}\, Q(u,u)| \leq \|\Delta u \|_{L^2(\Bbb{R}^2)}.$$
However, the absolute value of the corresponding sesquilinear  form $Q(u,v)$ cannot be majorized by 
$C\|\Delta u\|_{L^2(\Bbb{R}^2)} \| \Delta v\|_{L^2(\Bbb{R}^2)}$. Indeed,  the opposite assertion would yield an upper estimate of $\|r^{-1} \partial u/\partial r \|_{L^2(\Bbb{R}^2)}$ by the norm of $\Delta u$ in $L^2(\Bbb{R}^2)$ 
which is  wrong for a function linear near the origin.  $\hskip10mm \square$

\bigskip

\noindent
{\bf Remark 3} Note that under the additional orthogonality assumption
\begin{equation}\label{nn}
\int_0^{2\pi} u(r, \varphi) d\varphi =0 \quad {\rm for} \,\,\, r>0
\end{equation}
the above proof of Lemma \ref{lem4} provides inequality (\ref{1m}) with the sharp constant factor $3/4$ in the right-hand side. Besides, (\ref{nn}) implies
$${\rm Re}\int_{\Bbb{R}^2}\bigl(x_1 u_{x_1} + x_2 u_{x_2}\bigr) \Delta\overline{ u} \, \frac{dx}{|x|^2}\leq 0.\qquad \qquad \square$$

\bigskip

Using (\ref{1m}), we establish a new weighted G{\aa}rding type inequality. 

\begin{theorem}\label{th3}
Let $u\in C_0^\infty(\Bbb{R}^2 \backslash \{0\})$. Then inequality $(\ref{t1})$ holds.
\end{theorem}

\noindent
{\bf Proof.} Clearly, the right-hand side in (\ref{t1}) is equal to
$${\rm Re} \int_{\Bbb{R}^2 } \Delta\overline{u}\, \Delta\, (u \log |x|^{-1}) \, dx$$
$$=\int_{\Bbb{R}^2 } |\Delta u|^2 \log |x|^{-1}dx + 2{\rm Re} \int_{\Bbb{R}^2 }\Delta\overline{u} \cdot \nabla u\cdot \nabla \log |x|^{-1}dx.$$
Combining this identity  with (\ref{1m}), we arrive at the inequality
\begin{equation}\label{t2}
\int_{\Bbb{R}^2 } |\Delta u|^2 \log (e^2\, |x|)^{-1} dx \leq {\rm Re} \int_{\Bbb{R}^2 }\Delta\overline{u} \, \Delta(u \, \log |x|^{-1}) \, dx.
\end{equation}
Note that
$$\int_{\Bbb{R}^2 } \Delta u\cdot\Delta\overline{u}\cdot  \log |x|^{-1}  = - \int_{\Bbb{R}^2 } \nabla u \cdot \nabla(\Delta\overline{u}\cdot  \log |x|^{-1} ) \, dx $$
$$= {\rm Re}\int_{\Bbb{R}^2}\sum _{j=1}^2\Bigl(\nabla\frac{\partial u}{\partial x_j}\cdot \frac{\partial}{\partial x_j}\nabla \overline{u}\cdot \log |x|^{-1}  + \nabla u\cdot \frac{\partial}{\partial x_j}\nabla \overline{u}\cdot \frac{\partial}{\partial x_j}\bigl(\log |x|^{-1} \bigr) \Bigr)dx,$$
which is equal to 
$$\int_{\Bbb{R}^2}\Biggl(\sum _{j,k=1}^2\Bigl |\frac{\partial^2 u}{\partial x_j\partial x_k} \Bigr |^2\log |x|^{-1}  + \frac{1}{2}\sum _{j=1}^2\frac{\partial}{\partial x_j}|\nabla u|^2 \cdot \frac{\partial}{\partial x_j}\Bigl(\log |x|^{-1} \Bigr) \Biggr) dx.$$
Integrating by parts in the second term, we see that it vanishes. Thus, we conclude that
$$\int_{\Bbb{R}^2 } |\Delta u|^2 \log |x|^{-1} dx = \int_{\Bbb{R}^2 } |\nabla_2 u|^2\log |x|^{-1} dx$$
which together with (\ref{1m}) and the obvious identity
$$\int_{\Bbb{R}^2 } |\Delta u|^2  dx = \int_{\Bbb{R}^2 } |\nabla_2 u|^2dx$$
completes the proof of  (\ref{t1}). 

\smallskip

In order to see that no constant less than $1$ is admissible in front of the integral in the right-hand side of (\ref{t1}), it suffices to put 
$$u(x) = e^{i\langle x,\xi \rangle} \, \eta(x)$$ 
with $\eta\in C_0^\infty(\Bbb{R}^2)\backslash\{0\}$ into (\ref{t1}) and take the limit as $|\xi|\to \infty$.
$\hskip10mm \square$

\bigskip

\noindent
{\bf Remark 4} If the condition $u=0$ near the origin in Theorem \ref{th3} is removed, the above proof gives the additional term
$$\pi\Bigl( |\nabla u(0)|^2 -2\, {\rm Re} \bigl(u(0)\, \Delta\overline{u}(0)\bigr) \Bigr)$$
in the right-hand side of (\ref{t1}). $\hskip10mm \square$

\section{Dilation invariant Hardy's inequalities
 with \\
 remainder term}

\begin{theorem}\label{prop1}
$(i)$ Let $q$ denote a locally integrable nonnegative function on $(0,1)$. The best constant in the inequality
\begin{equation}\label{1u}
\int_{\Bbb{R}^n_+} x_n\, |\nabla u|^2 dx \geq C \int_{\Bbb{R}^n_+}q\Bigl(\frac{x_n}{(x_{n-1}^2 + x_n^2)^{1/2}}\Bigr)
\frac{|u|^2}{(x_{n-1}^2 + x_n^2)^{1/2}} dx, 
\end{equation}
for all $u\in C_0^\infty(\overline{\Bbb{R}^n_+})$,  which is equivalent to $(\ref{2u})$, is given by
\begin{equation}\label{3u}
\lambda:= \inf \frac{\displaystyle{\int_0^{\pi/2} \Bigl( \bigl | y'(\varphi)\bigr |^2 +\frac{1}{4}\bigl |y(\varphi)\bigr |^2\Bigr) \sin \varphi \, d\varphi}}{\displaystyle{\int_0^{\pi/2}\bigl |y(\varphi)\bigr |^2 q(\sin \varphi)\, d\varphi}},
\end{equation}
where the infimum is taken over all smooth functions on $[0,\pi/2]$.

$(ii)$ Inequalities $(\ref{1u})$ and $(\ref{2u})$ with a positive $C$ hold if and only if 
\begin{equation}\label{4u}
\sup\limits_{t\in (0,1)} (1- \log t) \int_0^t q(\tau)\, d\tau <\infty.
\end{equation}
Moreover,
\begin{equation}\label{40u}
\lambda \sim \Bigl (\sup\limits_{t\in (0,1)} (1- \log t) \int_0^t q(\tau)\, d\tau\Bigr)^{-1}, 
\end{equation}
where  $a\sim b$ means that $c_1 a\leq b\leq c_2 a$ with absolute positive constants $c_1$ and $c_2$. 
\end{theorem}

\smallskip

\noindent
{\bf Proof}  $(i)$ Let $U\in C_0^\infty(\overline{\Bbb{R}^2_+})$, $\zeta\in C_0^\infty(\Bbb{R}^{n-2})$, $x'=(x_1, \ldots x_{n-2})$, and let $N=const>0$. Putting
$$u(x) = N^{(2-n)/2}\zeta(N^{-1}x')U(x_{n-1}, x_n)$$
into (\ref{1u}) and passing to the limit as $N\to \infty$, we see that
 (\ref{1u}) is equivalent to the inequality
\begin{equation}\label{4}
\int_{\Bbb{R}^2_+} x_2 \bigl( |U_{x_1}|^2 + |U_{x_2}|^2\bigr) dx_1 dx_2 \geq C \int_{\Bbb{R}^2_+} q\Bigl(\frac{x_2}{(x_{1}^2 + x_2^2)^{1/2}}\Bigr)\frac{|U|^2dx_1 dx_2}{(x_1^2 + x_2^2)^{1/2}} , 
\end{equation}
where $U\in C_0^\infty(\overline{\Bbb{R}^2_+})$. 
Let $(\rho, \varphi)$ be the polar coordinates of the point $(x_1, x_2)\in {\Bbb{R}^2_+}$. Then (\ref{4}) can be written as
$$\int_0^\infty \int_0^{\pi} \Bigl( |U_\rho|^2 + \rho^{-2} |U_\varphi|^2 \bigr) \sin \varphi \, d\varphi \, \rho^2 d\rho \geq C \int_0^\infty \int_0^{\pi}|U|^2 q(\sin \varphi)\, d\varphi d\rho.$$
By the substitution
$$U(\rho, \varphi) = \rho^{-1/2} v(\rho, \varphi)$$
the left-hand side becomes
\begin{equation}\label{4a}
\int_0^\infty \int_0^{\pi}\Bigl ( |\rho v_\rho|^2 + |v_\varphi|^2 + \frac{1}{4} |v|^2\Bigr ) \sin \varphi \, d\varphi \frac{d\rho}{\rho} - {\rm Re}\int_0^{\pi}\int_0^\infty \overline{ v} \,  v_\rho \, d\rho \,  \sin \varphi \, d\varphi .
\end{equation}
Since $v(0) =0$, the second term in (\ref{4a}) vanishes. Therefore, (\ref{4}) can be written in the form
\begin{equation}\label{5}
\int_0^\infty \int_0^{\pi}\Bigl ( |\rho v_\rho|^2 + |v_\varphi|^2 + \frac{1}{4} |v|^2\Bigr ) \sin \varphi \, d\varphi \frac{d\rho}{\rho} \geq C \int_0^\infty \int_0^{\pi}|v|^2 q(\sin \varphi)\, d\varphi \frac{d\rho}{\rho}.
\end{equation}
Now, the definition (\ref{3u}) of $\lambda$ shows that (\ref{5}) holds with $C = \lambda$. 

\smallskip

In order to show the optimality of this value of $C$, put $t=\log \rho$ and $v(\rho, \varphi)= w(t,\varphi)$. Then (\ref{5}) is equivalent to 
\begin{equation}\label{6t}
\int_{\Bbb{R}^1} \int_0^{\pi}\bigl(|w_t|^2 + |w_\varphi|^2 + \frac{1}{4}|w|^2 \bigr)\sin \varphi \, d\varphi \, dt \geq C \int_{\Bbb{R}^1} \int_0^{\pi}|w|^2 q(\sin \varphi)d\varphi \, dt.
\end{equation}
Applying the Fourier transform $w(t, \varphi) \to \hat{w}(s,\varphi)$, we obtain
\begin{equation}\label{6}
\int_{\Bbb{R}^1} \int_0^{\pi}\Bigl(|\hat{w}_\varphi|^2 + \Bigl(|s|^2 + \frac{1}{4}\Bigr) |\hat{w}|^2\Bigr) \sin \varphi \, d\varphi \, ds\geq C \int_{\Bbb{R}^1} \int_0^{\pi}|\hat{w}|^2 q(\sin \varphi) d\varphi \, ds.
\end{equation}
Putting here
$$\hat{w}(s,\varphi) = \varepsilon^{-1/2}\eta(s/\varepsilon) y(\varphi),$$
where $\eta\in C_0^\infty(\Bbb{R}^1)$, $\|\eta\|_{L^2(\Bbb{R}^1)} =1$, and $y$ is a  function on $C^\infty ([0,\pi])$, and passing to the limit as $\varepsilon\to 0$, we arrive at the estimate
\begin{equation}\label{7}
\int_0^{\pi}\Bigl( |y'(\varphi)|^2 + \frac{1}{4} |y(\varphi)|^2\Bigr) \sin \varphi \, d\varphi  \geq C \int_0^{\pi} |y(\varphi)|^2 q(\sin \varphi) d\varphi
\end{equation}
where $\pi$ can be changed for  $\pi/2$ by symmetry. 
This together with (\ref{3u}) implies $\Lambda\leq \lambda$. 
 The proof of $(i)$ is complete.  
 
 \smallskip
 
 $(ii)$ Introducing the new variable $\xi = \log\cot \frac{\varphi}{2}$, we write (\ref{3u}) as
  \begin{equation}\label{5u}
\lambda = \inf_z \frac{\displaystyle{\int_0^\infty\Bigl( |z'(\xi)|^2 + \frac{|z(\xi)|^2} {4\, (\cosh\xi)^2}\Bigr) d\xi}}{\displaystyle{\int_0^\infty|z(\xi)|^2 \, q\Bigl( \frac{1}{\cosh \xi}\Bigr) \frac{d\xi}{\cosh \xi}}}. 
\end{equation}
Since 
$$|z(0)|^2 \leq 2 \int_0^1 \bigl(|z'(\xi)|^2 + |z(\xi)|^2\bigr) d\xi$$
and 
$$\int_0^\infty|z(\xi)|^2  \frac{ e^{2\xi}}{(1+ e^{2\xi})^2} d\xi \leq 2\int_0^\infty |z(\xi) - z(0)|^2 \frac{d\xi}{\xi^2} + 2 \, |z(0)|^2 \int_0^\infty \frac{ e^{2\xi}}{(1+ e^{2\xi})^2} d\xi$$
$$\leq 8\int_0^\infty |z'(\xi)|^2 d\xi + |z(0)|^2,$$
it follows from (\ref{5u}) that
\begin{equation}\label{1e}
\lambda \sim  \inf_z \frac{\displaystyle{\int_0^\infty |z'(\xi)|^2 d\xi + |z(0)|^2}}{\displaystyle{\int_0^\infty|z(\xi)|^2 \, q\Bigl( \frac{1}{\cosh \xi}\Bigr) \frac{d\xi}{\cosh \xi}}}\, .
\end{equation}
Setting $z(\xi)=1$ and $ z(\xi) = \min \{ \eta^{-1}\xi, 1\}$ for all positive $\xi$ and fixed $\eta>0$ into the ratio of quadratic forms in (\ref{1e}), we deduce that
$$\lambda \leq \min\Bigl\{ \Bigl(\int_0^\infty q\Bigl( \frac{1}{\cosh \xi}\Bigr) \frac{d\xi}{\cosh \xi} \Bigr)^{-1}, \,\, \Bigl(\sup\limits_{\eta>0} \eta \int_\eta^\infty q\Bigl( \frac{1}{\cosh \xi}\Bigr) \frac{d\xi}{\cosh \xi}\Bigr)^{-1} \Bigr \}.$$
Hence,
$$\lambda \leq c  \Bigl (\sup\limits_{t\in (0,1)} (1- \log t) \int_0^t q(\tau)\, d\tau\Bigr)^{-1}.$$

\smallskip

In order to obtain the converse estimate, note that
$$\int_0^\infty|z(\xi)|^2 \, q\Bigl( \frac{1}{\cosh \xi}\Bigr) \frac{d\xi}{\cosh \xi} $$
$$\leq 2\, |z(0)|^2 \int_0^\infty q\Bigl( \frac{1}{\cosh \xi}\Bigr) \frac{d\xi}{\cosh \xi}   + 2\int_0^\infty|z(\xi) - z(0)|^2 \, q\Bigl( \frac{1}{\cosh \xi}\Bigr) \frac{d\xi}{\cosh \xi} \, .$$

\noindent
The second term in the right-hand side is dominated by
$$ 8\sup\limits_{\eta >0} \Bigl( \eta \int_\eta^\infty q\Bigl( \frac{1}{\cosh \xi}\Bigr) \frac{d\xi}{\cosh \xi} 
\Bigr) \int_0^\infty |z'(\xi)|^2 d\xi$$
(see, for example, \cite{[M5]}, Sect. 1.3.1). Therefore,
$$\int_0^\infty|z(\xi)|^2 \, q\Bigl( \frac{1}{\cosh \xi}\Bigr) \frac{d\xi}{\cosh \xi}
 \leq 8\max \Biggl \{\int_0^\infty \!\! q\Bigl( \frac{1}{\cosh \xi}\Bigr) \frac{d\xi}{\cosh \xi}\, ,$$
$$  \sup\limits_{\eta >0}  \eta \! \int_\eta^\infty \!\! q\Bigl( \frac{1}{\cosh \sigma}\Bigr) \frac{d\sigma}{\cosh \sigma}\Biggr \} \!\Bigl(\int_0^\infty \! |z'(\xi)|^2 d\xi + |z(0)|^2 \Bigr)$$

\noindent
which together with (\ref{1e}) leads to the lower estimate
$$\lambda \geq \min\Bigl\{ \Bigl(\int_0^\infty q\Bigl( \frac{1}{\cosh \xi}\Bigr) \frac{d\xi}{\cosh \xi} \Bigr)^{-1}, \,\, \Bigl(\sup\limits_{\eta>0} \eta \int_\eta^\infty q\Bigl( \frac{1}{\cosh \xi}\Bigr) \frac{d\xi}{\cosh \xi}\Bigr)^{-1} \Bigr \}.$$
Hence,
$$\lambda \geq c  \Bigl (\sup\limits_{t\in (0,1)} (1- \log t) \int_0^t q(\tau)\, d\tau\Bigr)^{-1}.$$
The proof of $(ii)$ is complete. $\hskip10mm \square$

\medskip

Since (\ref{4u}) holds for $q(t) = t^{-1} (1- \log t)^{-2}$, Theorem \ref{prop1} (ii) leads to the following assertion.

\begin{corollary}\label{cor1}
There exists an absolute constant $C>0$ such that the  inequality
\begin{equation}\label{8u}
\int_{\Bbb{R}^n_+}  |\nabla v|^2 dx - \frac{1}{4} \int_{\Bbb{R}^n_+} \frac{|v|^2 dx}{x_n^2} \geq  
C \int_{\Bbb{R}^n_+}\frac{|v|^2 dx}{x_n^2\Bigl(1-\log \displaystyle{\frac{x_n}{(x_{n-1}^2 + x_n^2)^{1/2}}}\Bigr)^2}
\end{equation}
holds for all $v\in C_0^\infty (\Bbb{R}^n_+)$. The best value of $C$ is equal to
\begin{equation}\label{30}
\lambda:= \inf \frac{\displaystyle{\int_0^{\pi} \Bigl[ \bigl | y'(\varphi)\bigr |^2 +\frac{1}{4}\bigl |y(\varphi)\bigr |^2\Bigr] \sin \varphi \, d\varphi}}{\displaystyle{\int_0^{\pi}\bigl |y(\varphi)\bigr |^2(\sin \varphi)^{-1}\bigl(1-\log \sin \varphi)^{-2} d\varphi}},
\end{equation}
where the infimum is taken over all smooth functions on $[0,\pi/2]$. By numerical approximation,  $\lambda = 0.16 \dots$
\end{corollary}

A particular case of Theorem \ref{prop1} corresponding to $q =1$ is the following assertion.

\begin{corollary}\label{cor2}
The sharp value of $\Lambda$ in $(\ref{1})$ and $(\ref{2})$ is equal to
\begin{equation}\label{3}
\lambda:= \inf \frac{\displaystyle{\int_0^{\pi} \Bigl[ \bigl | y'(\varphi)\bigr |^2 +\frac{1}{4}\bigl |y(\varphi)\bigr |^2\Bigr] \sin \varphi \, d\varphi}}{\displaystyle{\int_0^{\pi}\bigl |y(\varphi)\bigr |^2 d\varphi}} \, ,
\end{equation}
where the infimum is taken over all smooth  functions on $[0,\pi]$. By numerical approximation, 
 $\lambda = 0.1564\ldots$
\end{corollary}

\noindent
{\bf Remark 5 } 
Let us consider the Friedrichs  extension $\tilde{{\cal L}}$ of the 
operator
\begin{equation}\label{3a}
{\cal L}: z \to -\bigl((\sin \varphi) z'\bigr)' + \frac{\sin\varphi}{4}z
\end{equation}
defined on smooth functions on $[0,\pi]$. It is a simple exercise to show that the energy space of $\tilde{{\cal L}}$ is compactly imbedded into $L^2(0, \pi)$. Hence, the spectrum of $\tilde{{\cal L}}$ is discrete and $\lambda$ defined by $(\ref{3})$ is the smallest eigenvalue of $\tilde{{\cal L}}$.

\noindent
{\bf Remark 6}
 The  argument used in the proof of Theorem \ref{prop1}  $(i)$ with obvious changes enables one  to obtain the following more general fact.
Let $P$ and $Q$ be measurable nonnegative functions in $\Bbb{R}^n$, positive homogeneous of degrees $2\mu$ and $2\mu -2$, respectively. The sharp value of $C$ in
\begin{equation}\label{6a}
\int_{\Bbb{R}^n} P(x)|\nabla u |^2 dx \geq C \int_{\Bbb{R}^n} Q(x)|u|^2 dx, \quad u\in C_0^\infty(\Bbb{R}^n),
\end{equation}
is equal to 
$$\lambda: =\inf\frac{\displaystyle{\int_{S^{n-1}} P(\omega) \Bigl( |\nabla_\omega Y|^2 + \bigl(\mu -1 +\frac{n}{2}\bigr)^2 |Y|^2\Bigr) }ds_\omega}{\displaystyle{\int_{S^{n-1}}  Q(\omega) |Y|^2 }ds_\omega}\, ,$$
where the infimum is taken  over all smooth functions on the unit sphere $S^{n-1}$. $\hskip 6mm \square$

\smallskip

A direct consequence of this assertion is the following particular case of  (\ref{6a}).

\medskip

\noindent
{\bf Remark 7}
 Let $p$ and $q$ stand for locally integrable  nonnegative functions on $(0,1]$ and let $\mu\in \Bbb{R}^1$.   If $n>2$, the best value of  $C$ in 
\begin{equation}\label{7a}
\int_{\Bbb{R}^n} |x|^{2\mu}p\bigl(\frac{x_n}{|x|}\bigr) |\nabla u |^2 dx \geq C \int_{\Bbb{R}^n} |x|^{2\mu -2} q\bigl(\frac{x_n}{|x|}\bigr)|u|^2 dx,
\end{equation}
where $u\in C_0^\infty(\Bbb{R}^n)$,  is equal to
\begin{equation}\label{7b}
\inf\frac{\displaystyle{\int_0^\pi\Bigl( |y'(\theta)|^2 + \bigl(\mu -1 +\frac{n}{2}\bigr)^2 |y(\theta)|^2\Bigr) p(\cos \theta) (\sin\theta)^{n-2} }d\theta}{\displaystyle{\int_0^\pi |y(\theta)|^2q(\cos \theta)(\sin\theta)^{n-2} }d\theta},
\end{equation}
with the infimum  taken over all smooth functions on the interval $[0,\pi]$. 

\smallskip

Formula (\ref{7b}) enables one to obtain a necessary and sufficient condition for the existence of a positive $C$ in (\ref{7a}). Let us assume that the function
$$\theta \to \frac{(\sin \theta )^{2-n}}{p(\cos \theta)}$$
is locally integrable on $(0,\pi)$. We make the change of variable $\xi = \xi(\theta)$, where
$$\xi(\theta) = \int_{\pi/2}^\theta \frac{(\sin \tau)^{2-n}} {p(\cos \tau)} d\tau, $$
and suppose that $\xi(0) = -\infty$ and $\xi(\pi) = \infty$. Then (\ref{7b}) can be written in the form
$$\lambda = \inf _z \frac{\displaystyle{\int_{\Bbb{R}^1} |z'(\xi)|^2 d\xi + \bigl(\mu -1+\frac{n}{2}\bigr)^2\int_{\Bbb{R}^1} |z(\xi)|^2\bigl(p(\cos \theta(\xi))\,  (\sin \theta(\xi))^{n-2} \bigr)^2 d\xi}}{\displaystyle{\int_{\Bbb{R}^1} |z(\xi)|^2  \, p(\cos \theta(\xi))\,  q(\cos \theta(\xi)) (\sin \theta(\xi))^{2(n-2)}\, d\xi}},$$
where $\theta(\xi)$ is the inverse function of $\xi(\theta)$. By Theorem $1$  in \cite{[M10]}, 
\begin{equation}\label{7ab}
\lambda \sim \inf\limits_{\xi\in \Bbb{R}^1\atop{ d>0, \delta>0}}\frac{\displaystyle{\frac{1}{\delta}} + \displaystyle{\int_{\theta(\xi -d -\delta)}^{\theta(\xi +d +\delta)} p(\cos \theta)\, (\sin \theta)^{n-2} \, d\theta}}{\displaystyle{\int_{\theta(\xi -d)}^{\theta(\xi +d )} q(\cos \theta)\, (\sin \theta)^{n-2} \, d\theta}}.
\end{equation}
Here the equivalence $a\sim b$ means that $c_1 b\leq a \leq c_2 b$, where $c_1$ and $c_2$ are positive constants depending only on $\mu$ and $n$. Hence (\ref{7a})  holds with a positive $\Lambda$ if and only if the infimum $(\ref{7ab})$ is positive. $\hskip10mm \square$

\smallskip

\noindent
{\bf Remark 8}
 In the case $n=2$ the best constant in $(\ref{7a})$ is equal to
\begin{equation}\label{7c}
\lambda: =\inf\limits_y \frac{\displaystyle{\int_0^{2\pi}\Bigl( |y'(\varphi)|^2 +\mu^2 |y(\varphi)|^2\Bigr) p(\sin \varphi) }d\varphi}{\displaystyle{\int_0^{2\pi} |y(\varphi)|^2q(\sin \varphi) }d\varphi},
\end{equation}
where the infimum is taken over all smooth functions on the interval $[0,2\pi]$. Note that (\ref{3u}) is a particular case of (\ref{7c}) with $\mu = 1/2$ and $p(t)=t$. 
As another application of (\ref{7c}) we obtain  the following special case of inequality (\ref{7a}) with $n=2$. 

\smallskip

{\it For all $u\in C_0^\infty(\Bbb{R}^2)$ the inequality}
\begin{equation}\label{7f}
\int_{\Bbb{R}^2}\frac{x_2^2\,  |u(x)|^2\, dx}{(x_1^2 +x_2^2)\displaystyle{\Bigl(\frac{\pi^2}{4} -\Bigl({\rm arcsin}\frac{x_1}{(x_1^2 + x_2^2)^{1/2}}\Bigr)^2\Bigr)}} \leq \frac{1}{2}\int_{\Bbb{R}^2} x_2^2\,  |\nabla u(x)|^2 dx
\end{equation}
{\it holds, where $1/2$ is the best constant.}

\smallskip

In order to prove (\ref{7f}),  we choose $p(t) = t^2$  and $ \mu =1$ in  (\ref{7c}), which  becomes
\begin{equation}\label{7d}
\lambda =\inf\limits_y\frac{\displaystyle{\int_0^{2\pi}\Bigl( |y'(\varphi)|^2 + |y(\varphi)|^2\Bigr) (\sin \varphi)^2 }d\varphi}{\displaystyle{\int_0^{2\pi} |y(\varphi)|^2q(\sin \varphi) }d\varphi},
\end{equation}
where $y$ is an arbitrary smooth $2\pi$-periodic function. Putting $\eta(\varphi) = y(\varphi)\sin\varphi$, we write (\ref{7d}) in the form
\begin{equation}\label{7e}
\lambda =\inf\limits_\eta \frac{\displaystyle{\int_0^{2\pi}|\eta'(\varphi)|^2 d\varphi}}{\displaystyle{\int_0^{2\pi}|\eta(\varphi)|^2 q(\sin\varphi)(\sin\varphi)^{-2}d\varphi}}
\end{equation}
with the infimum  taken over all $2\pi$-periodic functions satisfying $\eta(0) = \eta(\pi) =0$.  Let 
$$q(\sin\varphi) = \frac{(\sin\varphi)^2}{\varphi(\pi-\varphi)}.$$
In view of the well-known sharp inequality
$$\int_0^1\frac{|z(t)|^2}{t(1-t)} dt \leq \frac{1}{2}\int_0^1 |z'(t)|^2 dt$$
(see \cite{[HLP]}, Th. 262), we have $\lambda =2$ in (\ref{7e}). Therefore, (\ref{7a}) becomes (\ref{7f}). $\square$

\section{Generalized Hardy-Sobolev inequality with sharp constant}

Let ${\Omega}$ denote an open set in $\Bbb{R}^n$ and let $p\in [1, \infty)$. By $(p, a)$-capacity of a compact set $K\subset\Omega$ we mean the set function
$${\rm cap}_{p,a} (K, {\Omega}) = \inf\Bigl\{\int_{\Omega} |\nabla u|^p |x|^{-a} \, dx: \,\, u\in C_0^\infty({\Omega}), \,\, u\geq 1\  {\rm on}\, K\Bigr\}.$$
In the case $a=0$, $\Omega = \Bbb{R}^n$ we write simply ${\rm cap}_{p} (K)$.

\smallskip

The following inequality is a particular case of a more general one obtained in \cite{[M8]}, where $\Omega$ is an open subset of an arbitrary Riemannian manifold and $|\Phi(x, \nabla u(x)|$ plays the role of $|\nabla u(x)|\, |x|^{-a/p}$.

\begin{theorem}\label{th1}
{\rm (see $\cite{[M3]}$ for $q=p$ and $\cite{[M8]}$ for $q\geq p$)}

$(i)$ Let $q\geq p\geq 1$ and let ${\Omega}$ be an open set in $\Bbb{R}^n$. Then for an arbitrary  $u\in C_0^\infty({\Omega})$, 
\begin{equation}\label{60x}
\left(\int_0^\infty \bigl( {\rm cap}_{p, a} ({M_{t}}, {\Omega})\bigr)^{q/p}d(t^q)\right)^{1/q}\leq {\cal A}_{p,q}\left( 
\int_{\Omega}| \nabla u(x)|^p |x|^{-a} dx\right)^{1/p},
\end{equation}
where $M_t = \{x\in \Omega: \; |u(x)| \geq t\}$ and 
\begin{equation}\label{60y}
{\cal A}_{p,q}= \Bigl(\frac{\Gamma\bigl(\frac{pq}{q-p}\bigr)}{\Gamma\bigl(\frac{q}{q-p}\bigr)\Gamma\bigl(p\frac{q-1}{q-p}\bigr)}\Bigr)^{1/p-1/q} 
\end{equation}
for $q>p$, and 
\begin{equation}\label{60z}
{\cal A}_{p,p}=p(p-1)^{(1-p)/p}.
\end{equation}

$(ii)$ The sharpness of this constant is checked  by a sequence of radial functions in $C_0^\infty(\Bbb{R}^n)$.  
Moreover, there exists a radial optimizer vanishing at infinity.
\end{theorem}

 Being combined with the isocapacitary inequality
 \begin{equation}\label{63}
 \mu(K)^{\gamma} \leq \Lambda_{p,\gamma} \, {\rm cap}_{p, a} (K, {\Omega})
 \end{equation}
 where $\mu$ is a Radon measure in $\Omega$,  (\ref{60x}) implies the  estimate
  \begin{equation}\label{64}
 \Bigl( 
\int_0^\infty \bigl(\mu(M_t)\bigr)^{\gamma q/p} d(t^q) \Bigr)^{1/q} \leq {\cal A}_{p,q} \Lambda_{p,\gamma}^{1/p}\left( 
\int_{\Omega}| \nabla u(x)|^p|x|^{-a}dx\right)^{1/p}
\end{equation}
for all $u\in C_0^\infty(\Omega)$.

\smallskip

This estimate of $u$ in the Lorentz space ${\cal L}_{p/\gamma,q}(\mu)$ becomes the estimate in $L_q(\mu)$ for $\gamma = p/q$:
$$\|u\|_{L_q(\mu)} \leq {\cal A}_{p,q} \Lambda_{p,\gamma}^{1/p}\left( 
\int_{\Omega}| \nabla u(x)|^p|x|^{-a}dx\right)^{1/p}.
$$

\smallskip

In the next assertion we find the best value of $\Lambda_{p,\gamma}$ in (\ref{63}) for the measure $\mu = \mu_b$ defined by (\ref{65j}).

\begin{lemma}\label{lem33}
Let
\begin{equation}\label{69m}
1\leq p <n, \quad 0\leq a<n-p, \quad {\rm and} \,\,\,\,\,  a+p \geq b\geq \frac{an}{n-p}.
\end{equation}
 Then
\begin{equation}\label{39w}
\left(\int_{\Bbb{R}^n}\frac{dx}{|x|^b} \right)^{\frac{n-p-a}{n-b}}\leq \Bigl(\frac{p-1}{n-p-a}\Bigr)^{p-1} \frac{|S^{n-1}|^{\frac{b-p-a}{n-b}}}{(n-b)^{\frac{n-p-a}{n-b}}} \, {\rm cap}_{p,a}(K).
\end{equation}
The value of the constant factor in front of the capacity is sharp and the equality in $(\ref{39w})$ is attained at any ball centered at the origin.
\end{lemma}

\noindent
 {\bf Proof.} Introducing spherical coordinates $(r,\omega)$ with $r>0$ and $\omega\in S^{n-1}$, we have
 \begin{equation}\label{48j}
  {\rm cap}_{p, a}(K)= \inf\limits_{u\bigr |_K \geq 1}\int_{S^{n-1}}\int_0^\infty \Bigl( \Bigl |\frac{\partial u}{\partial r}\Bigr |^2 +\frac{1}{r^2}|\nabla_\omega u|^2\Bigr)^{\frac{p}{2}} r^{n-1-a}dr ds_\omega.
 \end{equation}
 Let us put here $r= \rho^{1/\varkappa}$, where
 $$\varkappa = \frac{n-p-a}{n-p}$$
  and $y = (\rho, \omega)$. The mapping $(r,\omega) \to (\rho, \omega)$ will be denoted by $\sigma$. Then (\ref{48j}) takes the form
 \begin{equation}\label{49j}
 {\rm cap}_{p, a} (K)= \varkappa^{p-1} \inf\limits_{v}  \int_{\Bbb{R}^n}\Bigl (\Bigl|\frac{\partial u}{\partial \rho}\Bigr |^2 + (\varkappa\rho)^{-2}|\nabla_\omega u|^2\Bigr)^{\frac{p}{2}} dy,
  \end{equation}
  where the infimum is taken over all $v= u\circ \sigma^{-1}$. Since $0\leq \varkappa\leq 1$ owing to the conditions $p<n$, $0<a<n-p$, and $a\geq 0$,  inequality (\ref{49j}) implies 
\begin{equation}\label{49v} 
  {\rm cap}_{p, a} (K)\geq \varkappa^{p-1} \inf\limits_{v}  \int_{\Bbb{R}^n} |\nabla u|^p dy \geq \varkappa^{p-1} {\rm cap}_{p} (\sigma(K))
 \end{equation} 
  which together with  the isocapacitary property of ${\rm cap}_{p} $ (see Corollary 2.2.3/2 \cite{[M5]}) leads to the estimate
\begin{equation}\label{50j} 
{\rm cap}_{p} (\sigma(K))\geq \Bigl(\frac{n-p}{p-1}\Bigr)^{p-1} |S^{n-1}|^{\frac{p}{n}} n^{\frac{n-p}{n}} \bigl( {\rm mes}_n(\sigma(K)\bigr)^{\frac{n-p}{n}}. 
  \end{equation}
Clearly,
$$\mu_b(K) = \frac{1}{\varkappa}\int_{\sigma(K)}\frac{dy}{|y|^\alpha}$$  
with 
\begin{equation}\label{50a} 
\alpha = n-\frac{n-b}{\varkappa} = \frac{b(n-p) -an}{n-p-a}\geq 0.
\end{equation} 
Furthermore, one can easily check that 
\begin{equation}\label{50b} 
 \mu_b(K)  \leq \frac{n^{1-\frac{\alpha}{n}}}{n-b} \, |S^{n-1}|^{\frac{\alpha}{n}}  \bigl( {\rm mes}_n(\sigma(K)\bigr)^{1-\frac{\alpha}{n}}
\end{equation} 
(see, for instance, \cite{[M8]}, Example 2.2). Combining (\ref{50b}) with (\ref{50j}), we find
\begin{equation}\label{50e}
\bigl( \mu_b(K)\bigr)^{\frac{n-p-a}{n-b}} \leq  \Bigl(\frac{p-1}{n-p}\Bigr)^{p-1}  \frac{|S^{n-1}|^{\frac{b-p-a}{n-b}}}{(n-b)^{\frac{n-p-a}{n-b}}} \, {\rm cap}_{p} (\sigma(K))
\end{equation} 
which together with (\ref{49v}) completes the proof of (\ref{39w}). $\hskip10mm \square$

\medskip

The main result of this section is as follows.

\begin{theorem}\label{th6}
Let conditions $(\ref{69m})$ hold and let $q\geq p$. 
Then for all $u\in C_0^\infty(\Bbb{R}^n)$
\begin{equation}\label{60a}
\Bigl(\int_0^\infty \bigl(\mu_b (M_t)\bigr)^{\frac{(n-p-a)q}{(n-b)p}} d(t^q)\Bigr)^{\frac{1}{q}} \leq {\cal C}_{p,q,a,b} \Bigl(\int_{\Bbb{R}^n}{| \nabla u(x)|^p}\frac{dx}{|x|^a} \Bigr)^{\frac{1}{p}},
\end{equation} 
where
\begin{equation}\label{60b}
{\cal C}_{p,q,a,b} = \Bigl(\frac{\Gamma\bigl(\frac{pq}{q-p}\bigr)}{\Gamma\bigl(\frac{q}{q-p}\bigr)\Gamma\bigl(p\frac{q-1}{q-p}\bigr)}\Bigr)^{\frac{1}{p}-\frac{1}{q}} \Bigl(\frac{p-1}{n-p-a}\Bigr)^{1-\frac{1}{p}}\Bigl(\frac{\Gamma(\frac{n}{2})}{2\pi^{\frac{n}{2}}(n-b)^{\frac{n-p-a}{p+a-b}}}\Bigr)^{\frac{p+a-b}{(n-b)p}}.
\end{equation} 
The constant $(\ref{60b})$   is best possible which can be shown by constructing a radial optimizing sequence in $C_0^\infty(\Bbb{R}^n)$. 
\end{theorem}

\noindent
{\bf Proof.} Inequality (\ref{60a})  is obtained by substitution of (\ref{60y})   and (\ref{39w}) into (\ref{64}). The sharpness of (\ref{60b}) follows from the part $(ii)$ of Theorem \ref{th1} and from the fact that the isocapacitary inequality (\ref{39w}) becomes equality for balls. $\hskip10mm \square$

\medskip

The last theorem contains  the best constant in the Il'in inequality (\ref{NN}) as a particular case 
 $q= (n-b)p/(n-p-a)$. We formulate this as the following assertion.
 
 \begin{corollary}\label{cor10}
 Let conditions $(\ref{69m})$ hold. Then for all $u\in C_0^\infty(\Bbb{R}^n)$
 \begin{equation}\label{60c}
\left(\int_{\Bbb{R}^n}{|u(x)|^{\frac{(n-b)p}{(n-p-a)}}}\frac{dx}{|x|^b} \right)^{\frac{n-p-a}{n-b}}\leq {\cal C}_{p,a,b} \Bigl(\int_{\Bbb{R}^n}{| \nabla u(x)|^p}\frac{dx}{|x|^a} \Bigr)^{\frac{1}{p}},
\end{equation}
where
$${\cal C}_{p,a,b} \!=\!\Bigl(\frac{p-1}{n\!-\!p\!-\!a}\Bigr)^{1-\frac{1}{p}}\Bigl(\frac{\Gamma(\frac{n}{2})}{2\pi^{\frac{n}{2}}(n-b)^{\frac{n-p-a}{p+a-b}}}\Bigr)^{\frac{p+a-b}{(n-b)p}}\Bigl(\frac{\Gamma\bigl(\frac{p(n-b)}{p-b}\bigr)}{\Gamma\bigl(\frac{n-b}{p-b}\bigr)\Gamma\bigl(1+\frac{(n-b)(p-1)}{p-b}\bigr)}\Bigr)^{\frac{p-b}{p(n-b)}}\!.$$
This constant     is best possible which can be shown by constructing a radial optimizing sequence in $C_0^\infty(\Bbb{R}^n)$.
\end{corollary}

\noindent

\section{Hardy's inequality with sharp Sobolev remainder term}

\begin{theorem}\label{prop4}
For all $u\in C^\infty(\overline{\Bbb{R}^n_+})$,  $u=0$ on $\Bbb{R}^{n-1}$, the sharp inequality
 \begin{equation}\label{7x}
\int_{\Bbb{R}^n_+} |\nabla u|^2 dx \geq \frac{1}{4}\int_{\Bbb{R}^n_+} \frac{|u|^2}{x_n^2}dx +\frac{\pi^{n/(n+1)}(n^2 -1)}{4\bigl(\Gamma\bigl(\frac{n}{2} +1\bigr)\bigr)^{2/(n+1)}} \|x_n^{-1/(n+1)} u\, \|^2_{L_{\frac{2(n+1)}{n-1}}(\Bbb{R}^n_+)}
\end{equation}
 holds. 
 \end{theorem}

\noindent
{\bf Proof.} We start with the Sobolev inequality
 \begin{equation}\label{8x}
 \int_{\Bbb{R}^{n+1}} |\nabla w|^2 dz \geq {\cal S}_{n+1} \|w\|^2_{L_{\frac{2(n+1)}{n-1}}(\Bbb{R}^{n+1})}
\end{equation}
with the best constant
 \begin{equation}\label{9x}
 {\cal S}_{n+1} = \frac{\pi^{(n+2)/(n+1)}(n^2 -1)}{4^{n/(n+1)}\bigl(\Gamma\bigl(\frac{n}{2} +1\bigr)\bigr)^{2/(n+1)}} 
\end{equation}
(see Rosen \cite{[R]}, Aubin \cite{[Au]}, and Talenti \cite{[Ta]}). 

\smallskip

Let us introduce the cylindrical coordinates $(r,\varphi, x')$, where $r\geq 0$, $\varphi\in [0, 2\pi)$, and $x'\in \Bbb{R}^{n-1}$. Assuming that $w$ does not depend on $\varphi$, we write (\ref{8x}) in the form
$$2\pi\int_{\Bbb{R}^{n-1}}\int_0^\infty \Bigl( \Bigl |\frac{\partial w}{\partial r} \Bigr |^2 + |\nabla_{x'} w|^2\Bigr) r\, dr dx' $$
$$\geq (2\pi)^{(n-1)/(n+1)} {\cal S}_{n+1} \Bigl(\int_{\Bbb{R}^{n-1}}\int_0^\infty |w|^{2(n+1)/(n-1)} r\, dr dx'\Bigr)^{(n-1)/(n+1)}.$$
Replacing $r$ by $x_n$, we obtain
$$\int_{\Bbb{R}^n_+} |\nabla w|^2 x_n\, dx \geq (2\pi)^{-2/(n+1)} {\cal S}_{n+1} \Bigl(\int_{\Bbb{R}^n_+} |w|^{2(n+1)/(n-1)} x_n\,  dx\Bigr)^{(n-1)/(n+1)}.$$
It remains to set here $w= x_n^{1/2}v$ and to use (\ref{9x}). $\square$

\end{document}